\documentclass{article}
\catcode`\@=11
\newtheorem{theorem}{Theorem}
\newtheorem{lemma}{Lemma}

\newtheorem{example}{Example}
\newtheorem{remark}{Remark}

\def\demo{\noindent{\bf Proof .-}}
\def\section{\@startsection {section}{1}{\z@}{-3.5ex plus -1ex
minus-.2ex}{2.3ex plus .2ex}{\normalsize\bf}}
\def\br{\hbox{\it I\hskip -2pt R}}

\def\bd{\hbox{\it I\hskip -2pt D}}
\def\bz{\hbox{\it Z\hskip -4pt Z}}

\def\bq{\hbox{\it C\hskip -5pt Q}}
\def\bc{\hbox{\it I\hskip -6pt C}}
\begin{document} 
\noindent
\begin{center}
{\Large\bf \textsc{On empty lattice simplices in dimension 4}}

\end{center}
\vskip.5truecm
\begin{center}
{Margherita Barile\footnote{Dipartimento di Matematica, Universit\`{a} di Bari ``Aldo Moro", Via E. Orabona 4, 70125 Bari, Italy, e-mail: barile@dm.uniba.it}\\ Dominique Bernardi\footnote{Universit\'e Pierre et Marie Curie, Institut Math\'ematique de Jussieu, 175 Rue du Chevaleret, F-75013, Paris, France, e-mail: bernardi@math.jussieu.fr}\\ Alexander Borisov\footnote{Department of Mathematics,
University of Pittsburgh, 301 Thackeray Hall, Pittsburgh, PA 15260, U.S.A., e-mail: borisov@pitt.edu}\\ Jean-Michel Kantor\footnote{Universit\'e Paris Diderot, Institut Math\'ematique de Jussieu, 175 Rue du Chevaleret, F-75013, Paris, France,  e-mail: kantor@math.jussieu.fr}}
\end{center}
\vskip1truecm
\noindent
{\bf Abstract} We give an almost complete classification of empty lattice simplices in dimension 4 using the conjectural results of Mori-Morrison-Morrison, later proved by Sankaran and Bober. In particular, all of these simplices correspond to cyclic quotient singularities, and all but finitely many of them have width bounded by 2.
\vskip0.5truecm
\noindent
{\bf Acknowledgments:} The research of the first author has been co-financed by the Italian Ministry of Education, University and Research (PRIN ``Algebra Commutativa, Combinatoria e Computazionale". The research of the third author has been supported by NSA, grant H98230-08-1-0129.

\section{Introduction}
Lattice polytopes have been extensively studied in recent years.
Their study shows strong connections with algebraic geometry and with number theory.
For example, the classical  Flatness Theorem  of Kronecker concerns the  arithmetic width of  compact convex sets, in particular lattice  polytopes. It
has many applications to the study of diophantine equations.
The program of classification of general algebraic varieties of  Mori and others gave an important role to singularities which are associated  with empty simplices, i.e., lattice simplices with no lattice points other than their  vertices. In dimension 3 the complete classification of such simplices is known, see \cite{White} and \cite{MS}.
Here we study   the  case of dimension 4. We give an almost complete classification (complete up to a finite number) of empty lattice simplices in dimension 4, by showing that all of them are cyclic. By the general result of \cite{B0} and the classification theorem of Sankaran \cite{S} and Bober \cite{Bob} this implies that all but finitely many of them belong to the families first discovered empirically by Mori-Morrison-Morrison \cite{MMM}.
As an application, we show that all but finitely many of the empty lattice simplices in dimension 4 have width bounded by 2 (see also \cite{HZ}, \cite{K}, and \cite{W}). 

\section{Preliminaries}
Let $\bd$ be a $d$-lattice which is a subgroup of $\bq^d$,  and let $\sigma$ be a $d$-dimensional rational convex cone spanned by $d$ linearly independent elements of $\bd$ in $\br^d$. Let $\bd'$ be the lattice spanned by  the intersections of the edges of $\sigma$ with $\bd$. 
\newline Let $e_1, \dots, e_d$ be the canonical basis of $\br^d$. Then the quotient group $G=\bd/\bd'$ is isomorphic to the direct sum $\oplus_{i=1}^d\bz/m_i\bz$ for some positive integers $m_1,\dots, m_d$.
According to the construction described, e.g., in \cite{F}, Chapter 2, the affine toric variety $U_{\sigma}$ associated with $\sigma$ can be viewed as the ($\bq$-factorial toric) quotient singularity $\bc^d/G$, where the action of $G$ on $\bc^d$ is defined as follows: for all $r_1\,\dots, r_d\in \bz$, and all $(z_1,\dots, z_d)\in \bc^d$,  
$$(\bar r_1,\dots, \bar r_d)(z_1,\dots, z_d)=(\rho_1^{r_1}z_1,\dots, \rho_d^{r_d}z_d),$$
\noindent
and, for all $i=1,\dots, d$, $\rho_i$ is a primitive $m_i$th root of unity. 
This quotient singularity is called {\it cyclic} if $G$ is a cyclic group.
  \begin{example}\label{example}{\rm (Cyclic case) Suppose $\bd=\sum_{i=1}^d\bz e_i$ and let $\sigma$ be the 
$d$-dimensional rational convex cone spanned by $v_1=-a_1e_1,$ $\dots,$ $v_{d-1}=-a_{d-1}e_{d-1},$ $v_d=-a_1e_1+\dots-a_{d-1}e_{d-1}+a_de_d$, where $a_1,\dots, a_d$ are non-zero integers such that $\gcd(a_1,\dots, a_d)=1$. Then $\bd'$ is the lattice spanned by $e_1,\dots, e_{d-1}, v_d$. Let $\phi:\br^d\to\br^d$ be the $\bq$-linear mapping such that $\phi(e_i)=e_i$ for all $i=1,\dots, d-1$ and $\phi(e_d)=\frac1{a_d}(a_1e_1+\cdots+a_{d-1}e_{d-1}+e_d).$
 Then $\phi(v_d)=e_d,$ so that $\tilde{\bd}'=\phi(\bd')=\sum_{i=1}^d\bz e_i$ and $\tilde{\bd}=\phi(\bd)=\tilde{\bd}'+\frac1{a_d}\bz(v_1+\dots+v_{d-1}-e_d)$. 
Hence $G$ is cyclic of order $a_d$, and the affine toric variety $U_{\sigma}$ associated with $\sigma$ can be viewed as the quotient singularity $\bc^d/G$, where $G$ is the group of complex $a_d$th roots of unity and its action  on $\bc^d$ is defined as follows: for all $\rho\in G$, 
$$\rho(z_1,\dots, z_d)=(\rho^{a_1}z_1,\dots, \rho^{a_d}z_d),$$
\noindent
where $\rho^{a_d}=1$. 
According to \cite{R}, this  quotient singularity is {\it terminal} if and only if the simplex 
$$\tilde{\Delta}=\{(x_1,\dots, x_d)\in\br^d\vert 0\leq x_i,\mbox{ for all }i=1,\dots,n,\sum_{i=1}^d x_i\leq 1\}$$
\noindent meets $\tilde{\bd}$ only in its vertices. But $\tilde{\Delta}=\phi(\Delta)$, where $\Delta$ is the simplex spanned by $e_1,\dots, e_{d-1}, v_d$. Hence the previous condition is equivalent to requiring that $\Delta$ is an empty lattice simplex.
}
\end{example}

\section{Cyclicity}
In this section we show that, in dimension 4, every terminal quotient singularity is cyclic.\newline
We first prove some results on good representatives of subspaces of finite vector spaces.\newline
Let $p$ be a positive prime number. For every $x\in  \bz_p$, let $s(x)\in\{0,1,\dots, p-1\}$ be such that $x$ is the residue class of  $s(x)$ modulo $p$. For all ${\bf x}=(x_1, x_2, \dots, x_n)\in  \bz_p^n$ we also set $s({\bf x})=\sum_{i=1}^ns(x_i).$ For every linear subspace $F$ of $\bz_p^n$, we consider the number
$$m(F)=\min_{{\bf x}\in F\setminus\{0\}}s({\bf x}).$$
White's Lemma (see, e.g., Corollary 1.4 in \cite{MS}) can then be formulated as follows.
\begin{lemma}\label{White}
Let $L$ be a line of $\bz_p^3$  and let $(a,b,c)$ be a direction vector of $L$. Then
$$m(L)\left\{\begin{array}{ll}
=p+1&\quad\mbox{if }abc\ne0\mbox{ and }(a+b)(a+c)(b+c)=0;\\
\leq p&\quad\mbox{otherwise}.
\end{array}
\right.
$$
\end{lemma}
From this lemma we can deduce a similar result in dimension two:
\begin{lemma}\label{White2}
Let $L$ be  a line of $\bz_p^2$, (where $p\ne 2$) and let $(a,b)$ be a direction vector of $L$. Then
$$m(L)\left\{\begin{array}{ll}
=p&\quad\mbox{if }a+b=0;\\
=\frac{p+1}2&\quad\mbox{if }(a+2b)(b+2a)=0;\\
\leq \frac{p-1}2&\quad\mbox{otherwise}.
\end{array}
\right.
$$
\end{lemma}
\demo We apply Lemma \ref{White} to the vector $(a,a+b, b)$. If $(a+b)(a+2b)(2a+b)\ne0$, we have $m(\bz_p(a,a+b,b))\leq p$. Hence there is some non-zero element $x\in\bz_p$ such that $s(xa)+s(xb)+s(x(a+b))\leq p$. But then $s(xa)+s(xb)< p$, whence we deduce that $s(x(a+b))=s(xa)+s(xb)$, so that $s(xa)+s(xb)\leq\frac{p}2$.
\par\smallskip\noindent
We can now prove the statement corresponding to White's Lemma in dimension four.
\begin{lemma}\label{White4}
If $P$ is a plane in $\bz_p^4$, then $m(P)\leq p$.
\end{lemma}
\demo For all $i\in\{1,2,3,4\}$ let $H_i$ be the hyperplane of $\bz_p^4$ defined by $x_i=0$. First assume that $P$ is contained in one of these hyperplanes, say $P\subset H_1$. Let ${\bf u}=(0,a,b,c)$ and ${\bf v}=(0,d,e,f)$ be the elements of a basis of $P$. If $a+b+c=0$, then $m(P)\leq m(\bz_p{\bf u})\leq p$ by Lemma \ref{White}. Otherwise, if we set
$${\bf w}=(a+b+c){\bf v}-(d+e+f){\bf u},$$
we have $m(P)\leq m(\bz_p{\bf w})\leq p.$ We now suppose that all intersections $L_i=P\cap H_i$ are lines.   Possibly rearranging the coordinates, we may assume that $L_1$ is spanned by ${\bf u}=(0,1,a,c)$. Then $L_2$ is spanned by ${\bf v}=(b,0,e,f),$ for some $a,b,c,e,f$ in $\bz _p.$ Then by application of Lemma \ref{White} to these lines, we may assume that $c=-a$ and one of $b$ and $e$ is equal to 1, so that $b=1$ and $e+f=0$, or $e=1$ and $b+f=0$.
Then ${\bf u}=(0,1,a,-a)$ and ${\bf v}=(1,0,b,-b)$ or ${\bf v}=(b,0,1,-b)$, where $ab\ne0$. If ${\bf v}=(1,0,b,-b)$, then ${\bf w}=b{\bf u}-a{\bf v}=(-a, b, 0,0)$ satisfies $s({\bf w})\leq p$ or $s(-{\bf w})\leq p$. Now suppose that ${\bf v}=(b,0,1,-b)$. The vector 
$${\bf w}={\bf u}-a{\bf v}=(-ab, 1,0, a(b-1))$$
is a direction vector of $L_3$. According to Lemma \ref{White}, we may assume that $b\ne 1$ and $a\ne 1$ and that either $ab=1$ or $a(1-b)=1$. The vector
$${\bf x}=b{\bf u}-a{\bf v}=(-ab, b,a(b-1),0)$$
is a direction vector of $L_4$. According to Lemma \ref{White}, we may assume that $b\ne a$ and $ab+b-a=0$. We deduce that $ab=1$, whence $a^2-a-1=0$ and $b=a-1$. In the latter case (which is impossible if $p\equiv \pm2$ (mod 5)), we thus have ${\bf v}=(a-1, 0,1,1-a)$. We set ${\bf y}=(1-a){\bf u}+(1+a){\bf v}=(a,1-a, a, 1-a)$. We then apply Lemma \ref{White2} to the line generated by $(a, 1-a)$ in $\bz_p^2$. Since $a\ne -1$ and $a\ne 2$, it contains a non-zero vector $(h,k)$ such that $s(h)+s(k)\leq \frac{p-1}2$ and then the vector ${\bf t}=(h,k,h,k)$ belongs to $P$ and we have $s({\bf t})\leq p$.
\par\smallskip\noindent
We can now prove the following
\begin{theorem} Every ${\bq}$-factorial toric terminal $4$-dimensional singularity is a cyclic quotient.
\end{theorem}
\demo
Suppose we have a ${\bq}$-factorial toric terminal $4$-dimensional singularity which is not a cyclic quotient. As above, denote by $G$ the quotient group of the ``big'' lattice $\bd$ by the ``small'' lattice $\bd'$. By the classification of the finite abelian groups, there is a prime $p$ such that $G$ contains a subgroup $G_1$ which is isomorphic to ${\bz}_p \times {\bz}_p$ (see \cite{H}, Theorem 3.3.3). The intermediate lattice that corresponds to this subgroup again defines a ${\bq}$-factorial toric terminal $4$-dimensional singularity. Thus we can assume that $G$ is isomorphic to ${\bz}_p \times {\bz}_p$. If we identify the lattice $\bd'$ with $\bz ^4,$ $G$ is a subgroup of $\frac{1}{p}\bz ^4 / \bz^4$. Multiplying the coordinates of all points by $p,$ we get a subgroup $P$ of $({\bz}_p)  ^4$ isomorphic to  ${\bz}_p \times {\bz}_p$. By Lemma \ref{White4} we can get a point  $A= (x_1,x_2,x_3,x_4)\in \bz ^4$ such that all $x_i \geq 0,$ $x_1+x_2+x_3+x_4\leq p,$  and modulo $p$ $A$ is a non-zero element of $P$. Hence, the point $\frac{1}{p}A\in \bd$ belongs to the simplex $\tilde{\Delta}$. This contradiction completes the proof.
\begin{remark}{\rm The result does not extend to higher dimensions. We give a counterexample in dimension 5. Let $\tilde{\bd}'=\bz^5$ and $\tilde{\bd}=\bz^5+\bz\frac1p(1,-1, 0,0,a)+\bz\frac1p(0,0,1,-1, b)$, where $a$ and $b$ are integers not divisible by $p$. Then $\tilde{\bd}/\tilde{\bd}'\simeq \bz_p\times \bz_p$ is not cyclic. We show that it gives rise to a terminal singularity. The lattice points of $\tilde{\bd}$ are those of the form
$$Q=(x_1, x_2, x_3, x_4, x_5)=(\alpha_1+\frac{\beta}p, \alpha_2-\frac{\beta}p, \alpha_3+\frac{\gamma}p, \alpha_4-\frac{\gamma}p, \alpha_5+\frac{a\beta+b\gamma}p),$$
where $\alpha_i\in\bz$ for $i=1,\dots, 5$ and $\beta,\gamma\in\bz$. We show that if $Q\in\tilde{\Delta}$, then $Q\in\tilde{\bd}'$. Set $x=\sum_{i=1}^5x_i=\sum_{i=1}^5\alpha_i+\frac{a\beta+b\gamma}p$, and suppose that $x_i\geq 0$ for all $i=1,\dots, 5$.   Note that $x_1\in \bz$ iff $x_2\in \bz$ and $x_3\in \bz$ iff $x_4\in \bz$. If $x_i\not\in \bz$ for $i=1,\dots,  4$, then $\beta$ and $\gamma$ are non-zero, so that one of $\alpha_1,\alpha_2$ and one of $\alpha_3,\alpha_4$ must be positive, that is, at least 1. Then  $x\geq2$, so that $Q\not\in \tilde\Delta$. If $x_i\in \bz$ for $i=1,2$ and $x_i\not\in \bz$ for $i=3,4$, 
then $p$ divides $\beta$ but not $\gamma$, so that $\gamma\neq0$ and, once again, $\alpha_3+\alpha_4\geq 1$. Moreover, $p$ does not divide $a\beta+b\gamma$, so that $x_5>0$. Hence $x>1$. Finally, suppose that $x_i\in \bz$ for $i=1,\dots, 4$. Then $p$ divides both $\beta$ and $\gamma$, and $x_5\in \bz$. Hence $Q\in \tilde{\bd}'$. 
  } 

\end{remark}
\section{Width}
In this section we prove the following result.
\begin{theorem} Up to possibly a finite number of exceptions, every empty simplex in dimension 4 has width 1 or 2.
\end{theorem}
\demo 
Like in Example 1 and the proof of Theorem 1, we identify the lattice $\bd '$ with $\bz ^4$ so that that our simplex is the standard simplex $\Delta$. Then the lattice $\bd$ is obtained from $\bd '$ by adding multiples of one rational point, with coordinates $\frac{1}{N} (a_1,a_2,a_3,a_4).$  Here $a_i$ are integers, and $N$ is a nonzero integer.

By the general result of \cite{B0} cyclic terminal quotients  form  a finite number of families (see \cite{B0} for a precise definition of a family). All positive-dimensional families have been found empirically in \cite{MMM}. They are those associated with the vectors $U$ that are  obtained by applying the $j$th projection to the integer vectors $U'$ of the following types:
\begin{list}{}{}
\item{(i)} $U'=(x, -x, y, z, -y-z)$, or
\item{(ii)} $U'=(x, -2x, y, -2y, x+y)$, or
\item{(iii)} belonging to the ones listed in Table 1.9 in \cite{MMM}, which we reproduce below.
\end{list}
$$
\begin{array}{lc}
\mbox{Stable Quintuple}&\mbox{Linear Relations}\\
(9,1,-2,-3,-5)&(0,2,1,0,0),(1,1,0,0,2),(2,0,1,2,2)\\
(9,2,-1,-4,-6)&(0,1,2,0,0),(0,2,0,1,0),(2,0,2,1,2)\\
(12,3,-4,-5,-6)&(0,2,0,0,1),(1,0,0,0,2),(1,2,2,2,0)\\
(12,2,-3,-4,-7)&(0,2,0,1,0),(1,1,0,0,2),(2,0,2,1,2)\\
(9,4,-2,-3,-8)&(0,1,2,0,0),(0,2,0,0,1),(2,0,2,2,1)\\
(12,1,-2,-3,-8)&(0,2,1,0,0),(1,2,0,2,1),(2,0,1,2,2)\\
(12,3,-1,-6,-8)&(0,2,0,1,0),(1,0,0,2,0),(1,2,2,0,2)\\
(15, 4, -5,-6,-8)&(0,2,0,0,1), (2,0,2,2,1)\\
(12,2,-1,-4,-9)&(0,1,2,0,0),(0,2,0,1,0),(2,0,2,1,2)\\
(10,6,-2,-5,-9)&(0,2,1,2,0),(1,0,0,2,0),(1,2,2,0,2)\\
(15, 1, -2,-5,-9)&(0,2,1,0,0),(2,0,1,2,2)\\
(12,5,-3,-4,-10)&(0,2,0,0,1),(0,2,2,1,0),(2,0,2,2,1)\\
(15, 2, -3,-4,-10)&(0,2,0,1,0),(2,0,2,1,2)\\
(6,4,3,-1,-12)&(0,2,2,2,1),(2,0,0,0,1)\\
(7,5,3,-1,-14)&(0,2,2,2,1),(2,0,0,0,1)\\
(9,7,1,-3,-14)&(0,2,0,0,1),(2,0,2,2,1)\\
(15,7,-3,-5,-14)&(0,2,0,0,1),(2,0,2,2,1)\\
(8,5,3,-1,-15)&(0,2,2,1,1),(2,0,0,1,1)\\
(10,6,1,-2,-15)&(0,0,2,1,0),(2,2,0,1,2)\\
(12,5,2,-4,-15)&(0,0,2,1,0),(2,2,0,1,2)\\
(9,6,4,-1,-18)&(0,2,2,2,1),(2,0,0,0,1)\\
(9,6,5,-2,-18)&(0,2,2,2,1),(2,0,0,0,1)\\
(12,9,1,-4,-18)&(0,2,0,0,1),(2,0,2,2,1)\\
(10,7,4,-1,-20)&(0,2,2,2,1),(2,0,0,0,1)\\
(10,8,3,-1,-20)&(0,2,2,2,1),(2,0,0,0,1)\\
(10,9,4,-3,-20)&(0,2,2,2,1),(2,0,0,0,1)\\
(12,10,1,-3,-20)&(0,2,0,0,1),(2,0,2,2,1)\\
(12,8,5,-1,-24)&(0,2,2,2,1),(2,0,0,0,1)\\
(15,10,6,-1,-30)&(0,2,2,2,1),(2,0,0,0,1)\\

\end{array}
$$
\noindent 
The linear relations on the right are the linear relations of the entries in the corresponding quintuple with coefficients 0,1 and 2. The original interest in these relations comes from algebraic geometry (cf. \cite{MMM}). The fact that they exist for each quintuple is of utmost importance to us. Note also that in each quintuple the sum of the entries is zero. For any integer $j$ between 1 and 5, the $j$th projection is the operation of omitting the $j$th coordinate. 

For example, the quintuple  $(9,1,-2,-3,-5)$ with the 5th projection produces the singularity generated by the extra point $\frac{k}{n}(9,1,-2,-3).$ One can actually choose $k$ to be $1$, and there are natural restrictions on $n$, but we will not need the detailed analysis of these families. The families $(i)$  and $(ii)$ give a singularity when one specifies $x$ and $y$ to be any rational numbers and chooses a coordinate to drop. Once again, there are certain restrictions on $x$ and $y$ for the lattice $\bd$ not to contain any points in $\Delta,$ but we will not need these restrictions.

 The completeness of this empirical classification of \cite{MMM} was first established by Sankaran in \cite{S}. A conceptually different proof was recently obtained by Bober \cite{Bob}. It also follows from a result of Vasyunin \cite{Vasyunin}. Because the total number of families is finite, all but finitely many terminal cyclic quotients are obtained in the way described above.

To prove that the corresponding simplex has width at most two, we will produce a linear map from $\bd$ to $\bz$ that sends vertices of $\Delta$ to the set $\{0,1,2\}$ or $\{-1,0,1\}$.  In all cases this map will have integer coefficients and will send the extra point to $0$. For example, for the singularity defined by adding the point $\frac{k}{n}(9,1,-2,-3),$ we define the map by $2x_2+x_3$. This can be traced to the first linear relation in the above table, given by the coefficients $(0,2,1,0,0)$. This works because we dropped the 5th coordinate, and the corresponding entry in the relation was $0$. Note that the vertex $(0,0,0,0)$ is sent to $0,$ together with the vertices $(1,0,0,0)$ and $(0,0,0,1)$. The vertex $(0,1,0,0)$ is sent to $2,$ and the vertex $(0,0,1,0)$ is sent to $1$. This will work for each quintuple when we drop a coordinate with entry $0$ for one of the relations in the table. Fixing any relation, if the entry is not $0$, but $2$, we can get another relation by subtracting the given one from $(2,2,2,2,2),$ and then do the same. Finally, if the entry is $1,$ we can subtract from the relation $(1,1,1,1,1)$ and get a relation with coefficients $-1,$ $0,$ and $1.$ The same construction then produces a map from $\bd$ to $\bz$ that maps the vertices of $\Delta$ to $\{-1,0,1\}$.

The same idea works for the cases $i$ and $ii$. In fact, for the case $i$ the width is $1$, using the relation $(1,1,0,0,0,)$ and for the case $ii$ we can use one of the relations $(2,1,0,0,0)$ and $(0,0,2,1,0)$. This completes the proof.

{\bf Remark 1} The above theorem cannot be extended to all 4-dimensional empty lattice simplices. Consider the simplex with vertices $(0,0,0,0)$,  $(1,0,0,0)$,  $(0,1,0,0)$, $(0,0,1,0)$, and $(6,14,17,65).$ It has width 4, and is the only simplex of width 4 among simplices of determinant at most $1000$ (cf. \cite{HZ}). 

{\bf Remark 2} Even though dropping one of the five coordinates from a given quintuple gives five different cyclic quotient singularities, the corresponding empty simplices are the same. Indeed, this simplex can be described as sitting in the affine subspace of $\br^5$ of points with sum of coordinates 1, with vertices $(1,0,0,0,0)$, $(0,1,0,0,0)$, $(0,0,1,0,0)$, $(0,0,0,1,0)$, and $(0,0,0,0,1)$. The lattice is the restriction to this affine subspace of a lattice in $\br ^5$ that is obtained from $\bz ^5$ by adding multiples of $\frac{k}{n}(a_1,a_2,a_3,a_4,a_5),$ where  $(a_1,a_2,a_3,a_4,a_5)$ is the given quintuple. When we drop a coordinate, we project to $\br ^4,$ sending one of the vertices of the simplex to $(0,0,0,0).$ This projection is an isomorphism between the lattice described above and the lattice described in the proof of Theorem 2.

 \end{document}